\documentclass[11pt, a4paper]{article}

\usepackage{fullpage}
\usepackage{amsmath,amssymb,srcltx}
\usepackage{amsmath,amssymb,amsthm}
\usepackage{graphicx}
\usepackage{amscd}
\usepackage{amsmath}
\usepackage{amsfonts}
\usepackage{amssymb}
\usepackage{setspace}
\usepackage{enumerate}         % better lists
\usepackage{color}
\usepackage{url}
\usepackage{amsthm}
\usepackage{hyperref}
\usepackage{amsmath}
\usepackage{bm}
\usepackage{xy}
\usepackage{stmaryrd}
\usepackage{fancyhdr}
\usepackage{authblk}

\theoremstyle{plain}
\newtheorem{theorem}{Theorem}[section]

\newtheorem{lemma}[theorem]{Lemma}

\newtheorem{proposition}[theorem]{Proposition}

\newtheorem{definition}[theorem]{Definition}

\newtheorem{assumption}[theorem]{Assumption}
\theoremstyle{remark}
\newtheorem{remark}[theorem]{Remark}

\numberwithin{equation}{section}

\newcommand{\RR}{\mathbb{R}}

\newcommand{\NN}{\mathbb{N}}

\newcommand{\lzp}{{\bf{L}}^0_+}

\newcommand{\fps}{(\Omega,\cF,(\cF_t)_{\tT},\massP)}

\newcommand{\LL}{\mathbf{L}}
\newcommand{\massP}{\mathbf{P}}
\newcommand{\massQ}{\mathbf{Q}}
\newcommand{\massE}{\mathbb{E}}

\newcommand{\la}{\lambda}

\newcommand{\cC}{\mathcal{C}}

\newcommand{\cF}{\mathcal{F}}

\renewcommand{\i}{\infty}

\newcommand{\tT}{0 \leq t \leq T}

\newcommand{\dom}{\operatorname{dom}}

\providecommand{\keywords}[1]{\textbf{{Keywords:}} #1}
\providecommand{\subjclass}[1]{\textbf{{MSC2010:}} #1}

\begin{document}

\title{Non-concave utility maximization problem with transaction costs and a given consistent price system \thanks{The authors gratefully acknowledge financial support from National natural science foundation of China, youth science fund project under grant 11801365 and 11901097, from Educational research projects for young and middle-aged teachers in Fujian under grant JT180073.} }
\author[1]{Lingqi Gu}
\author[2]{Yiqing Lin\thanks{Corresponding author: yiqing.lin@sjtu.edu.cn}}
%\author[3]{Junjian Yang}
\affil[1]{\small College of Mathematics and Statistics, Fujian Normal University, 350117 Fuzhou, China}
\affil[2]{\small School of Mathematical Sciences, MOE-LSC, Shanghai Jiao Tong University, 200240 Shanghai, China}
%\affil[3]{\small Fakult\"at f\"ur Mathematik und Geoinformation, TU Wien, A-1040 Wien, Austria}

\date{\small \today}

\maketitle

\begin{abstract}\noindent
We investigate expected utility maximization problems from the terminal liquidation value in continuous time in markets with  transaction costs and one fixed consistent price system, where a non-concave utility function is defined on the positive half real line.
The  sufficient  conditions are given by  the convex conjugate  of the  utility function, then   the existence of the optimizer is proved by a maximizing sequence. Finally, we show that the value function of the envelope of the utility function and the concave envelope of the value function are  coincide.
\end{abstract}

\noindent \keywords{Non-concave utility maximization problem, Transaction costs, Consistence price system, Asymptotic elasticity   }

  \vspace{3mm}

\noindent \subjclass{91B16, 91G70,93E20, }

\maketitle

\section{Introduction}
In this paper, we study the existence of optimal portfolio for maximizing expected utility of the liquadation value at the terminal time in markets with proportional transaction costs $\lambda$ and one given fixed  $\lambda-$local constistent price 
system.
Preferences of an investor are discribed by non-concave function defined on the positive real axis, then the existence is proved by reducing the original problem to the traditional concave case with the help of the concave envelope of the utility function.
  
For the classical concave case, the existence of the optimizer is well known.
There are three different useful methods to solve classical utility maximization problems: Martingale methods or duality approach (see for example \cite{CK96, KS99}), dynamic programming (see for instance \cite{RS05}) and BSDE (see for example\cite{HIM05}).
 In markets with transaction costs, it is difficult to use dynamic programming if we consider the expected utility of the terminal liquidation value, hence duality approach is applied more in this setting, such as \cite{CMKS14,CS16duality,CSY17,Sch17,CPSY18}.   While the  problems constructed by Kabanov model is solved by martingale methods such as  in \cite{Kab99,KS09, CO11}. 

The case that the utility function is  not concave attracted attention after Kahneman and Tversky \cite{KT79}  proposed prospect theory, in which the utility function is  illustrated  as S-shaped based on the empirical research, as they
believed  risk aversion and risk appetite can be mutually translated.
 Then, Berkelaar,Kouwenberg and Post \cite{BKP04} studied the utility optimization problem under S-shaped utility in continuous time.  
 Similarly in continuous time setting, Carassus and Pham \cite{CP09} investigated such problem but assuming the utility function is piecewise concave.
 Reichlin\cite{Rei13} relaxed the restriction of the shape of the utility function, assumed the utility function is increasing, upper semi-continuous satisfying the growth  condition. He proved the existence of the optimal strategies by maximizing sequences. 
All of them have the assumption that the market is complete.
 %As our main reference, we will borrow some idea of \cite{Rei13} to solve our problem.
In the incomplete market with non-concave utility on the positive half real line, due to the mathematical difficulties, the research restricts itself in  discrete-time  setting. 
In \cite{CRR15}, a sufficient condition including asymptotic elasticity on  non-concave, non-differentiable and random utility function is provided. 
Blanchard, Carassus and R{\'a}sonyi \cite{BCR18} applied  dynamic programming with measurable selection arguments to establish the no-arbitrage condition and proved the existence of the optimal strategies in the similar setting of   \cite{CRR15}. There are many other branches (behavior economics, incentive schemes, aspiration utility and goal reaching problems) of literature that investigate non-concave problems which will not be listed here.   

As we know until now, all existing research results on optimal investment under frictions are obtained under the assumption of  a concave utility function, only \cite{PPR17} involves  non-concave utilities in market with friction, in which transaction costs are assumed to be superlinear functions of the trading speed. 
We are working in the setting of \cite{CSY17} but remove the  concavity and differentiability of the utility function. 
From the idea of \cite{Rei13}, we give a sufficient condition which grantees the existence of the optimal solution  by concave  envelope of the utility function and its convex conjugate.  
 
The paper is structured as follows. Financial model is introduced in section 2, where we give the precise definition of the non-concave utility function and the original utility maximization problem. In section 3, as a useful tool, the envelope of the utility function and the convex conjugate of utility function are detailed introduced. The main theorem of the existence of the optimizer is proved in section 4 and the properties of the value function are studied in Section 5. Finally, some useful techniques of convex analysis are introduced in appendix.

\section{Formulation of the utility maximization problem}
Fix a finite time horizon $T > 0$.
We consider a scalar-valued strictly positive, c\`adl\`ag, adapted price process  $S=(S_t)_{0 \leq t \leq T}$ based on a fixed filtered probability space $\fps$ satisfying the usual conditions.
Additionally, we assume $\mathcal{F}_{T-}=\mathcal{F}_T$ and $S_{T-} =S_T$.
    \footnote{These assumptions avoid special notations for possible trading at the terminal time $T$.
              In this case, we may assume without loss of generality that the agent liquidates her position in the stock shares at time $T$.
              For more details see e.g.~\cite[Remark 4.2]{CS06} or \cite[p.~1895]{CS16duality}.}
The market involves proportional transaction costs $0<\lambda< 1$, i.e., an investor buys stock shares at the higher ask price $S_t$ and only receives a lower bid price $(1-\la)S_t$ when selling them.

Fix $0<\lambda<1$. A {\it self-financing} trading strategy, modeling holdings in units of the bond and of stock, is an $\mathbb{R}^2$-valued, predictable process $\varphi = (\varphi^0_t,\varphi^1_t)_{0\leq t\leq T}$ of finite variation
   such that
     $$ \int_s^td\varphi_u^0 \leq -\int_s^tS_ud\varphi_u^{1,\uparrow}+ \int_s^t(1-\lambda)S_ud\varphi_u^{1,\downarrow}, \mbox{ for all } 0\leq s<t\leq T, $$
  where $\varphi^1=\varphi^{1,\uparrow}-\varphi^{1,\downarrow}$ denotes the canonical decompositions of $\varphi^1$ into the difference between two increasing processes.
 The {\it liquidation value} $V_t^{liq}(\varphi)$ of a trading strategy $\varphi$ at time $t\in [0, T]$ is defined by
    $$ V_t^{liq}(\varphi) := \varphi^0_t+(\varphi_t^1)^+(1-\la)S_t-(\varphi_t^1)^-S_t.$$
 A trading strategy $\varphi=(\varphi^0,\varphi^1)$ is {\it admissible}, if  the liquidation value $V_t^{liq}(\varphi)\geq 0 \ a.s.$,  for all $0\leq t\leq T$.
 For $x>0$, $\mathcal{A}(x)$ denotes the set of all admissible, self-financing trading strategies $\varphi$ starting with the initial endowment $(\varphi_{0}^0,\varphi_{0}^1)=(x,0)$,
   and $\mathcal{C}(x)$ denotes the convex subset of terminal liquidation values as
  \begin{equation*}
   \mathcal{C}(x):= \left\{V^{liq}_T(\varphi):\,\,\varphi\in\mathcal{A}(x)\right\} = \left\{\varphi_T^0:\,\,\varphi=(\varphi^0,\varphi^1)\in\mathcal{A}(x),\ \varphi_T^1=0 \right\} \subseteq L_+^0(\massP).
  \end{equation*}
\begin{remark}
 Actually, by super-hedging theorem in \cite{sch14super}, the set $\mathcal C(x)$ can be  described by $\{V^{liq}_T\in \mathbf L_+^0: \mathbb E^{\massQ} [ V^{liq}_T ] \leq x  \}$.
\end{remark}
\vspace{3mm}
As usual,  the utility function will be used to describe the risk appetite of investors,
but the difference here is that a investor  can change her perception of risk, so we relax the usual restriction of concavity of the utility function, even the differentiability.

\begin{definition}
A utility is a function $U:(0,\infty)\to \mathbb R$ with $U(\infty)>0$, which is non-constant,increasing,upper semi-continuous satisfies the growth condition
\begin{equation}\label{growthcond}
\lim_{x\to \infty}\frac{U(x)}{x}=0.
\end{equation}

\end{definition}
  We also consider utility functions defined on the positive praxis, that is, taking the value $-\infty $ on $(-\infty,0)$.
We define $U(0):=\lim_{x\searrow 0}U(x)$ and  $U(\infty):=\lim_{x\nearrow \infty}U(x)$.
In the concave case, the growth condition \eqref{growthcond} is equivalent to the Inada condition at $\infty$ that $U'(\infty)=0$.
We assume that $U(\infty)>0$, because adding a constant to $U$ does not change the preferences described by $U$.
Moreover, the condition \eqref{growthcond} together with the assumption $U(\infty)>0$ imply that there is always a concave function $g:\mathbb R \to \mathbb R \cup \{-\infty\}$
satisfying $g \geq U$.

Given the initial endowment $x>0$, the investor wants to maximize her {\it expected utility} at terminal time $T$:
 \begin{equation}\label{J5}
   \massE\big[U\big(V_T^{liq}(\varphi)\big)\big] \to \max!, \qquad \varphi\in\mathcal{A}(x).
 \end{equation}
 In this case, the original problem can also be rewritten as 
 \begin{equation}\label{originalprob}
 u(x, U)=\sup\{E[U(f):f\in \mathcal C(x)]\}.
 \end{equation}
We define $E[U(f)]:=-\infty$ if $U^{-}(f)\notin \mathbf L^1$.
\begin{definition} \label{CPS}
Fix $0<\la <1$.
A \textnormal{$\la$-consistent price system} is a two-dimensional strictly positive process $Z=(Z^0_t,Z^1_t)_{0\le t\le T}$
 with $Z^0_0=1$, consisting of a martingale $Z^0$ and a local martingale $Z^1$ under $\massP$, such that
  \begin{equation*}
    \widetilde{S}_t:=\frac{Z^1_t}{Z^0_t} \in [(1-\la)S_t, S_t],\quad {\rm a.s.},\quad 0\leq t\leq T.
  \end{equation*}
%Denote by $\mathcal{Z}^{\lambda} (S)$ the set of $\la$-consistent price systems, and we say that $S$ satisfies the condition $(CPS^\la)$ of {admitting a $\la$-consistent price system} if $\mathcal{Z}^{\lambda}(S)$ is nonempty .

     The identification is given by the   $Z_T^0=\frac{d\massQ}{d\massP}$, where the probability measure $\massQ$ on $\cF $ is equivalent to $\massP$.

  %  We say that $S$ satisfies the condition $(CPS^\la)$ \textnormal{locally}
%if there exists a strictly positive process $Z$ and a sequence $\{\tau_n\}_{n\in\NN}$ of $[0,T] \cup\{\i\}$-valued stopping times, increasing to infinity,
%such that each stopped process $Z^{\tau_n}$ defines a $\lambda$-consistent price system for the stopped process $S^{\tau_n}$.
%This process $Z$ is called a local $\lambda$-consistent price system, and $\mathcal{Z}^{loc,\lambda}$ denotes the set of all such processes.
\end{definition}
%\begin{problem}
%should we assume that the fixed Z is a CPS or local-CPS? \\
%if CPS Z, the question is: $\mathcal B(x)$ is too small to have the real optimizer?becasue
%in \cite{CSY17}, a two way crossing is used and then get the main result:under the condition  $\inf_{0 \leq t \leq T}\widehat V_t^{liq}>0$ is the dual optimizer $\widehat Y \in \mathcal B(x)$ a local martingale,that means the first component is a local martingale.\\
%if use local-CPS, then can we use $Z_T^0=\frac{d\massQ}{d\massP}$ with a probability measure $\massQ \approx \massP$?
%\end{problem}
 %By super-hedging theorem in \cite{sch14super},we can define a a probability measure $\massQ \approx \massP$, such that for every local-CPS
Throughout the whole paper, we have the following assumption.
\begin{assumption}\label{fixedCPS}
 We consider a  {$\la$-consistent price system} $Z$ with  $Z_T^0=\frac{d\massQ}{d\massP}$, where $\massQ$ is  a probability measure equivalent to  $ \massP$.
\end{assumption}

\section{Existence of the optimizer}

One may notice that the sufficient condition called {\it{reasonable asymptotic elasticity}}, i.e., $AE(U)<1$, introduced in \cite{KS99} is not mentioned here. As the utility function is not concave, we need to find a similar condition to replace the condition  of $AE(U)<1$.
For this, we will apply some convex analysis technique, especially the concave envelope of the utility function.

%\section{Concave envelope of $U$}
 % We emphasize that  the following  useful lemmas are proved in \cite[section 2.2]{Rei13},we write down only for the convenience of the reader.
\begin{definition}
The concave envelope of the utility function $U$ is the smallest concave function 
$U_c:\mathbb R \to \mathbb R \cup \{-\infty\}$ such that $U_c(x)\geq U(x)$ holds for all $x \in \mathbb R.$
\end{definition}

With the definition of the concave envelope $U_c$ of $U$, we can conclude some properties of $U_c$
which  are proved in \cite[Appendix A]{Rei13}.

\begin{lemma}\cite[Lemma 2.8]{Rei13}
 The concave envelope $U_c$  of $U$ is finite, continuous on $(0,\infty)$ and satisfies the growth condition \eqref{growthcond}. 
The set $\{U < U_c \}:=\{x \in \mathbb R_+ | U(x)<U_c (x)\}$ is open and its countable
connected components are bounded open intervals. Moreover, $U_c$  is locally affine on the set $\{U < U_c \}$, in the sense that it is affine on each of the above intervals.
\end{lemma}

\begin{proof}
See \cite[Appendix]{Rei13}.
\end{proof}

The relationship between $U$ and $U_c$ can be studied by the key tool {\it the convex conjugate of $U$} defined by
\begin{equation}
V(y):=\sup_{x>0}\{U(x)-xy\}.
\end{equation}
\begin{remark}
If $U$ is  not concave,then the concave envelope $U_c$ is not strictly concave, which means   the function $V$ is no longer smooth; we therefore work with the subdifferential which is 
denoted by $\partial V$ for the convex function $V$ and by $\partial U_c$  for the concave function $U_c$. The right- and left-hand derivatives of $V$ are denoted by $V_+^{'}$   
and $V_-^{'}$.
\end{remark}

\begin{lemma}\cite[lemma 2.9]{Rei13}\label{chaV}
 The function $V$ is convex, decreasing, finite on $(0,\infty)$ and satisfies $V(x) = \infty$ for $x < 0$. The utility $U$ and its concave envelope $U_c$ have the same convex conjugate.
Moreover, it holds that
\begin{equation}\label{sameconvexconj}
U_c(x)-xy=V(y) \Longleftrightarrow x \in -\partial V(y) \Longleftrightarrow y \in \partial U_c(x)
\end{equation}
\end{lemma}
\begin{proof}
Monotonicity of $V$  as well as $ V >  - \infty$ on  $(0,\infty)$ and $V= \infty$ on
$(-\infty,0)$ follow from the definition. The fact that $V < \infty$ on  $(0,\infty)$ follows from the growth
condition \eqref{growthcond}. The equivalence follows from Lemma \ref{helplemma}(v) and (vi).
%\cite[Theorem 23.5]{Roc70}.  \cite[Theorem 11.1]{RW98}
Lemma \ref{helplemma} (ii) gives that $U$ and $U_c$ have the same convex conjugate. 
\end{proof}

     To make the  utility maximization problem solvable, one sufficient condition introduced in \cite{KS99} is (RAE) reasonable asymptotic elasticity of the concave utility function, i.e.,  $AE(U)=\limsup_{x \to \infty}\frac{xU'(x)}{(U)}<1$.
In our case, 
because of the non-concavity of the utility function,
 the condition of RAE should be redefined.
Fortunately, \cite[Lemma 6.3]{KS99} gives some equivalent inequality to the asymptotic elasticity $AE(U)$.
Then, \cite{DPT01} shows the other formulation constructed by the convex conjugate for the  asymptotic elasticity, which is used by Reichlin \cite{Rei13} to solve the  non-concave utility maximization problem in a complete market. 
We also apply this idea to solve our  non-concave utility maximization problem in  markets with  transaction costs.
We define an extending asymptotic elasticity by the conjugate $V$.
\begin{equation}
EAE(V):=\limsup_{y\to 0}\sup_{q \in \partial V(y)}\frac{|q|y}{V(y)}.
\end{equation}
 
\begin{lemma}\label{Vine}\cite[Lemma 2.10]{Rei13}
The asymptotic elasticity condition $EAE(V)<\infty$ is equivalent to the existence
of two constants $\gamma >0$ and $y_0>0$ such that 
$$V(\mu y)\leq \mu^{-\gamma}V(y),  \forall\ \mu \in (0,1] \mbox{and } y \in (0,y_0].$$
Moreover, if $EAE(V)<\infty$ is satisfied, 
then there are two constants $x_0>0 \mbox{ and }  k > 0$ such that 
\begin{equation}\label{bdUc}
0 \leq U_c(x)\leq kU(x) \mbox{ on }  (x_0, \infty).
\end{equation}
\end{lemma}  
\begin{proof}
The proof of the equivalence can be seen for example in \cite[lemma 4.1]{DPT01}.
For the second part, we refer to \cite[proof of Lemma 2.10]{Rei13}. The following proof is only for the convenience of the reader.
We first prove that there are $x_0>0 \mbox{ and } \gamma<1$ such that 
\begin{equation}\label{boundedUc}
0 \leq U_c(x)\leq U(x)+ \gamma U_c(x)  \mbox{  on   }  (x_0, \infty).
\end{equation}
By \cite[Proposition 4.1]{DPT01} 
we know $EAE(V)<\infty$, together with the growth condition \eqref{growthcond} show the existence of two constants  $x_0  \mbox{ and } \gamma<1$ such that $\sup_{y\in \partial U_c(x) \leq \gamma U_c(x)}$ holds on $(x_0, \infty)$.
By moving $x_0$ to the right  if necessary,
we may assume that $U(x_0)$ is positive.
Moreover,
$U_c$ is locally affine on $\{U< U_c\}=\cup_i(a_i,b_i)$.
Hence, for $x \in (a_i,b_i)\cap (x_0, \infty)$,
we rewrite $U_c(x)$,
use that $U_c(a_i)=U(a_i)\leq U(x)$ and apply for $y=U'_c(x)$ the above inequality $U'_c(x)x\leq \gamma U_c(x)$ to get
$$U_c(x)=U_c(a_i)+U'_c(x)(x-a_i)\leq U(x)+U'_c(x)x\leq U(x)+\gamma U_c(x).$$
Set $k=\frac{1}{1-\gamma}$,then $U_c(x) \leq kU(x)$.
For $x\in \{U=U_c\}\cap (x_0, \infty)$,
\eqref{boundedUc} follows since $U(x)=U_c(x)$ is positive.
\end{proof}
In the traditional case with concave function, i.e., $U=U_c$, the  assumption that $u(x,U)< \infty$   for some $x >0$ is given  to exclude the trivial case. But this is not  sufficient to guarantee the existence of a maximizer, 
even in  the concave case, see for example Section 5 in \cite{KS99} and \cite{BG11}, 
so a stronger assumption will be proposed here.

\begin{assumption}\label{vfinite}
 Suppose that $  \massE[V(yZ_T^0)] <\i, \forall \ y>0$.
\end{assumption}

\begin{lemma}{\cite[Lemma 3.2]{Rei13}}\label{UVequ}
Under the assumption $EAE(V)<\i$, Assumption \ref{vfinite} is equivalent to $u(x,U)<\i$,
for some $x>0$.
\end{lemma}
\begin{proof}
For the proof we refer to \cite[Appendix C]{Rei13}.
We write down only for the convenience of the reader.
Because  $\massE[U(f)]\leq \massE[U_c(f)]\leq \massE[V(yZ_T^0)]+yx_0 $, 
for all $f\in \mathcal C(x_0) \mbox{ and } y>0$,
Assumption \ref{vfinite} implies $u(x,U)<\i$ for some $x>0$.
For the other direction,
we  need to prove Assumption \ref{vfinite} 
 from  $AE(V)<\i$ and $u(x,U_c)<\i$ for some $x>0$(see for instance \cite[Lemma 5.4]{WZ09}) . 
%It is therefore sufficient to show that $EAE(V)<\i$ and $u(x,U)<\i$ for some $x>0$$
%imply $u(x,U_c)<\i$ for some $x>0$.
Let $x_0>0$ and $k$ be as  given in \eqref{bdUc} in Lemma \ref{Vine}, i.e., 
\begin{equation}
0 \leq U_c(x) \leq kU(x) \mbox{ on } (x_0, \infty). \label{UcleqkU}\end{equation} 

If $U >0$  on $(0, \infty)$, 
  fix some $f \in \mathcal C(x)$ and apply \eqref{UcleqkU} on the set $\{ f > x_0 \}$, which 
gives $$\massE[U_c( f )] \leq U_c(x_0 )+k\massE[U( f )].$$ 
Taking the supremum over all $f \in \mathcal C(x)$ implies $$u(x, U_c) \leq U_c(x_0 )+ku(x,U).$$
Hence,  if   $u(x, U )<\i$ for some $x>0$, then we have $u(x, U_c )< \i$.

If $U(0)<0$, choose $\epsilon$ small enough such that $x-\epsilon>0$, 
fix $f \in \mathcal C(x-\epsilon)$  and apply the above argument to $f_{\epsilon}:=f+\epsilon$ and $U(x)-U(\epsilon)$.
This gives 
\begin{align*}
\massE[U_c(f)]&\leq \massE[U_c(f_{\epsilon})]\\
	              &\leq U_c(x_0) +k \massE[(U(f_{\epsilon})-U({\epsilon}))\textbf{1}_{  \{f_{\epsilon} \geq {x_0}\} }]+k\massE[U(\epsilon)\textbf{1}_{  \{f_{\epsilon} \geq {x_0}\} }] \\
&\leq U_c(x_0) +k \massE[U(f_{\epsilon})-U({\epsilon})]+k\massE[U(\epsilon)\textbf{1}_{  \{f_{\epsilon} \geq {x_0}\} }] \\
&\leq U_c(x_0) +k u(x, U)+k|U(\epsilon)|,
\end{align*}
where $f_{\epsilon} \in \mathcal C(x )$ 
is used in the last step. Taking the supremum over all $f \in \mathcal C(x-\epsilon)$, we obtain $u(x-\epsilon, U_c) \leq ku(x,U)+\tilde k$ for some constant $\tilde k$ and the result follows.

\end{proof}

\begin{theorem}{(compare \cite[Theorem 3.3]{Rei13})}\label{exop}
Under Assumption \ref{vfinite},  
there exists some $\hat f \in \mathcal C(x_0)$ such that $u(x_0, U)=\massE[U(\hat f)]$, for all $   x_0 \in (0, \i)$,.
\end{theorem}
We also try to  work directly with a maximizing sequence $(f^n)_{n \in \mathbb N}$.  Koml\'os-type argument  do not work because of the lack of the concavity of the utility function $U$,  so we first proposed the following statement.
\begin{proposition}{(compare \cite[Proposition 3.4]{Rei13})}\label{helpTheo}
Under Assumption \ref{vfinite}, consider a   sequence $(f^n)$ with $f^n \in \mathcal C(x_n)$, where $x_n$ is a sequence satisfying $x_n \to x >0 .$
There is some $\hat f \in \mathcal C(x) $ such that $$\limsup_{n \to \i}\massE[U(f^n)]\leq E[U(\hat f)].$$
\end{proposition}

\begin{proof}
We first prove the uniformly integrability of the family $(U^+(f^n))_{n \in \mathbb N}$.
It is clear if $U$ is bounded from above. 
Hence, we may assume $U(\i)=\i$.
The sequence $(x_n )$ is bounded by $x_0$, and it follows that 
\begin{align*}
\massE[U^+(f^n)\textbf{1}_{  \{U^+(f^n) > \alpha\} }]&\leq \massE[(V^+(yZ_T^0)+yf^nZ_T^0)\textbf{1}_{  \{U^+(f^n) > \alpha\} }]\\
&\leq \massE[V^+(yZ_T^0)\textbf{1}_{  \{U^+(f^n) > \alpha\} }]+yx_0.
\end{align*}
Therefore, for any $y>0$, we have
$$\lim_{\alpha \to \i}\sup_{n}\massE[V^+(yZ_T^0)\textbf{1}_{  \{U^+(f^n) > \alpha\} }] =0. $$
Since  $V^+(yZ_T^0) \in \mathbf L^1$ is trivially uniformly integrable,
 we only need to show that $\sup_{n}\massP[U^+(f^n) > \alpha] \to 0 \mbox{ for  } \alpha \to \i$.
 For this, fixed a sequence $\alpha_i \to \i$ and let $\tilde x:=\inf\{x>0|U(x)>0\}$ denote the first point where the utility becomes positive.
By definition of $\tilde x$ and $x_0$,
we have $U^+(f^n) \leq U^+(\tilde x+f^n) $ and $\tilde x+f^n \in \mathcal C(\tilde x + x_0)$.
By the proof of Lemma \ref{UVequ}, we have $u(x, U_c)<\i$ for all $x>0$,
we get
$$\sup_n \massP[U^+(f^n) > \alpha_i]\leq \sup_n \frac{ \massE[U^+(f^n)]}{\alpha_i}\leq \frac{u(x_0+\tilde x,U)}{\alpha_i}  \longrightarrow 0,$$
which completes the proof of the first part.

By passing to a subsequence that realizes the  $lim sup$,
we may assume  the sequence $(\massE[U(f^n)])$ converges and its limit is denoted by,
  $\gamma_1$.
Again by passing to a further subsequence still relabelled as $f^n$  that realizes $\limsup \massE[Z_T^0f^n]$, we can also assume that $(\massE[Z_T^0f^n])$ converges to some limit $\gamma_2$, it is clear $\gamma_2 \leq x$.
From the proof of the first part, 
we know that the family $([U^+(f^n))_n$ is uniform integrable;
the same is trivially true for  $((Z_T^0f^n)^-)_n$.
By \cite[Corollary 3.9]{Bal84},
it follows then that there exist $\hat g_1 \in \mathbf{L}^1$ and $\hat g_2 \in \mathbf{L}^1$ such that 
\begin{equation}\label{L1bdg}
\massE[\hat g_1] \geq \gamma_1 \mbox{ and }   \massE[\hat g_2] \leq \gamma_2,  
\end{equation}
and for a.e. $\omega \in \Omega$,
 there exists a subsequence $n_k(\omega)$ such that 
\begin{equation}\label{limUfnk}
\lim_{k\to \i}U(f^{n_k(\omega)}(\omega))=\hat g_1(\omega) \mbox{ and } 
\lim_{k\to \i}Z_T^0(\omega)f^{n_k(\omega)}(\omega)=\hat g_2(\omega)
\end{equation}
since $Z_T^0>0 \ \massP-$a.s., we can define $\hat f(\omega):=\frac{\hat g_2(\omega)}{Z_T^0(\omega)}$, which shows $\hat f Z_T^0=\hat g_2$. Then, from \eqref{L1bdg}, 
it follows that $\massE[Z_T^0\hat f]=\massE[\hat g_2] \leq \gamma_2 \leq x$ 
which means $\hat f \in \mathcal C(x)$.
Moreover, it follows from \eqref{limUfnk} that
$$\lim_{k\to \i}f^{n_k(\omega)}(\omega)=\lim_{k\to \i}\frac{Z_T^0(\omega)f^{n_k(\omega)}(\omega)}{Z_T^0(\omega)}=\frac{\hat g_2(\omega)}{Z_T^0(\omega)}=\hat f(\omega).$$
Together with the upper-semicontinuity of $U$ and \eqref{limUfnk} we obtain 
$$U(\hat f(\omega))\geq \lim_{k\to \i}U(f^{n_k(\omega)}(\omega))=\hat g_1(\omega).$$
Taking expectations and using \eqref{L1bdg} gives $\massE[\hat f]\geq \massE[\hat g_1]\geq \gamma_1$. 
We conclude that $\hat f$ is a maximizer since $\hat f \in \mathcal C(x)$.
 \end{proof}

The proof of Theorem \ref{exop} is a direct application of Proposition \ref{helpTheo}.

\begin{proof}[Proof of Theorem \ref{exop}]
 Consider a maximizing sequence  $(f^n)_{n\in \NN} $ in $\mathcal C(x_0)$.
$U$ is increasing, so we can assume that the constraint $\massE[Z_T^0f^n]=x_0$ is satisfied for each $n \in \NN$. 
Proposition \ref{helpTheo} gives some $\hat f\in \mathcal C(x_0)$ such that
$\massE[U(\hat f)] \geq \lim_n\massE[U(f^n)]=u(x_0,U)$, 
which shows $\hat f$ is a maximizer.
\end{proof}
There are some good properties of the optimizer in the traditional case with concave utility such as the uniqueness, see for example \cite{KS99, CSY17}. Do these properties still stay in our case?
In the following,  
we will  study properties of the optimizer.

\begin{remark}[non-uniqueness]
After having the existence of the solution to $u(x, U)$,  a natural question on the uniqueness of the optimizer  arises. Actually, the optimizer is not necessarily unique.
Reichlin gave an example (see \cite[Example 3.6]{Rei13}) to show that even without transaction costs in a complete market, the solution of the non concave utility maximization problem could not be unique.
\end{remark}

\begin{lemma}{(compare \cite[Lemma 3.7]{Rei13})}
 If $U$ is continuously differentiable with $\{U <U_c\}=\cup_{i=1}^n(a_i,b_i)$  for some
fixed $n$, let $\hat f$ be a maximizer for $u(x, U)$,
then there is $ y> 0$ such that  $\hat f$ satisfies $U'(\hat f)=yZ_T^0$ on  $\{\hat f>0\}$.
\end{lemma}
\begin{proof}
Fix $\epsilon>0$, define $A_{\epsilon}:=\{\hat f>\epsilon\}$.  
Consider the following maximization problem
\begin{align} \label{maxEUf}
\max\massE[U(f)\mathbf 1_{A_{\epsilon}}], 
 \mbox{ subject to } & f\in \mathbf L_+^0,  \mbox{  and   }
\massE[fZ_T^0\mathbf 1_{A_{\epsilon}}] \leq \massE[\hat fZ_T^0\mathbf 1_{A_{\epsilon}}].
\end{align}
If there is some element $\tilde f$ in \eqref{maxEUf} with $\massE[U(\tilde f)\mathbf 1_{A_{\epsilon}}] > \massE[U(\hat f)\mathbf 1_{A_{\epsilon}}]$,
the candidate $f':=\tilde f\mathbf 1_{A_{\epsilon}}+\hat f\mathbf 1_{A_{\epsilon}^c}$ is feasible for the problem $u(x,U)$ and satisfies $\massE[U(f')]>\massE[U(\hat f)]$ which contradicts the optimality of $\hat f$.
Hence $\hat f$ solves \eqref{maxEUf}.

Now fix some $f \in \mathbf L^{\i}$ and define $f_y:=\hat f + y(f-c)\mathbf 1_{A_{\epsilon}}$
where $c$ is set by $c:=\frac{\massE[Z_T^0f\mathbf 1_{A_{\epsilon}}]}{\massE[Z_T^0\mathbf 1_{A_{\epsilon}}]}$.
Note that
 $$\massE[Z_T^0f_y\mathbf 1_{A_{\epsilon}}] = \massE[Z_T^0\hat f\mathbf 1_{A_{\epsilon}}]+y(\massE[Z_T^0f\mathbf 1_{A_{\epsilon}}]-c\massE[Z_T^0\mathbf 1_{A_{\epsilon}}])=\massE[Z_T^0\hat f\mathbf 1_{A_{\epsilon}}]$$
holds for every $y$.
Moreover $f\in \mathbf L^{\i}$ implies $f_y  \geq 0$ on $A_{\epsilon}$ for $y$ enough small.
Hence $f_y$ is a feasible candidate for the problem \eqref{maxEUf} and this yields
$$\limsup_{y \to 0}\frac{\massE[(U(  f_y)-U(\hat f))\mathbf 1_{A_{\epsilon}}]}{y}\leq 0.$$
Since $U$  is continuously differentiable and concave on $(b_n ,\i), U'$ is bounded on $(\epsilon,\i)$.
By the mean value theorem  $\frac{ (U(  f_y)-U(\hat f))\mathbf 1_{A_{\epsilon}}}{y}$  is bounded by a constant on $A_{\epsilon}$.
  Interchanging limit and expectation, then we obtain $$0 \geq \massE[U'(\hat f)(f-c)\mathbf 1_{A_{\epsilon}}].$$
Replacing $f$ by $- f$ shows that the expectation must vanish. 
Define  $\gamma:=\frac{\massE[U'(\hat f) \mathbf 1_{A_{\epsilon}}]}{\massE[Z_T^0\mathbf 1_{A_{\epsilon}}]}$, then
we see that $$\massE\left[\left(U'(\hat f) - Z_T^0 \gamma\right)f\mathbf 1_{A_{\epsilon}}\right]=0, \forall f\in \mathbf L^{\i},$$
which implies that $U'(\hat f)=Z_T^0\gamma$ on $A_{\epsilon}$.
The same approach for $\tilde \epsilon\in (0,\epsilon)$ gives $U'(\hat f)=Z_T^0 \tilde \gamma$ on $A_{\tilde \epsilon}$ for some constant $\tilde \gamma$.
Since $A_{  \epsilon}  \subset A_{\tilde \epsilon}$,
we have that $U'(\hat f)=Z_T^0\gamma=Z_T^0 \tilde \gamma$ on $A_{\epsilon}$ and we infer $\gamma =\tilde \gamma.$
This can be done for any $\tilde \epsilon>0$ and we obtain  $U'(\hat f)=Z_T^0\gamma$ on $\cup_{\epsilon>0}A_{\epsilon}=\{\hat f>0\}$, which completes the proof.
\end{proof}

\section{The value function}
 We may be curious about some properties of the   value function $u(x,U)$, especially if $u(x, U)$ has the same properties as the utility function $U$, i.e., is $u(x, U)$    increasing,upper-semicontinuous? 
On the other hand, we are interested in the concave envelope of $u(x,U)$ and their relationship.

\begin{proposition}{\cite[Proposition 4.2]{Rei13}}\label{upsemicont}
Under Assumption \ref{vfinite} the value function $u(x,U)$ is upper-semicontinuous. 
If $U$ is in addition continuous, then the value function $u(x,U)$ is continuous
on $(0,\i)$.
\end{proposition}
\begin{proof}
The proof please see also \cite[proof of Proposition 4.2]{Rei13}. 
For upper-semicontinuity, consider a sequence $x_n \searrow x \in (0, \i)$.
As  the existence of the maximizer is guaranteed by Theorem \ref{exop}, let  $\hat f^n$ be  the maximizer for  $u(x_n, U)$.  By Proposition \ref{helpTheo}, there is some $\hat f \in \mathcal C(x)$
such that 
$$\limsup_{n \to \i}u(x_n, U)=\limsup_{n\to \i}\massE[U(f^n)] \leq \massE[U(\hat f)] \leq u(x,U),$$   
which completes the proof of upper-semicontinuity for $x \in (0,\i)$.
Moreover, it is known from  \cite[Theorem 6]{Sio16}  that $u(x_n, U_c) \searrow U_c(0)$
for $x_n \searrow 0$ and this implies $u(x_n, U ) \searrow U(0)$ for $x_n \searrow 0$
since $U_c(0)=U(0) \leq u(x, U) \leq u(x, U_c)$ for all $x>0$.

We still need to show the lower-semicontinuity of $u(x,U)$  if $U$ is continuous.
Fix a sequence $x_n \nearrow x \in (0, \i)$  and a random variable  $f_0 \in \mathcal C(x_0)$ with $U^-(f_0) \in \mathbf{L}^1$. 
Then we construct  a sequence $(f_n)_n$ with $f_n \in \mathcal C(x_n), f_n \nearrow f_0$ and 
$U^-(f_n) \in \mathbf{L}^1, \mbox{ for all } n$. 
With the upper bound $U^-(f_1)$,  we can use Dominated convergence for $(U^-(f_n))_n$, and  
Fatou's lemma for $(U^+ ( f_n ))_n$,
then
obtain  $$ \massE[U(f_0)]\leq \liminf_{n \to \i}\massE[U(f_n)]\leq \liminf_{n\to \i}u(x_n, U).$$
Taking the supremum over $f_0 \in \mathcal C(x)$ finishes the proof.
\end{proof}

In the traditional case with concave utility, Duality Theorem  (see for example \cite{KS99, CSY17}) shows that $u$ and $v$ are conjugate. We wonder if $U$ is not concave, can this relationship stay.

\begin{theorem}{(compare \cite[Theorem 4.1]{Rei13})}
 For any utility function $U$ with convex conjugate $V$, if Assumption \ref{vfinite} is satisfied, the value function $u(x, U_c )$ and the concave envelope of $u(x, U)$ coincide on 
$(0, \i)$ and we have
\begin{equation}\label{defv}
v(y):=\massE[V(yZ_T^0)]=\sup_{x>0}\{u(x,U)-xy\}, y>0.
\end{equation}
\end{theorem} 
\begin{remark}
Actually, \eqref{defv}  tells us the relationship of $u $ and $v$ that  $v$ is the conjugate of the value function $u(x,U)$ even in the case of non-concave utility  in markets with transaction cost but under one fixed consistent price system, because in this setting, utility function do not need to apply the mimimax theorem to exchange the supremum over payoffs and the infimum over pricing densities. 
 
\end{remark}
\begin{proof}
To show \eqref{defv}, we first  claim that
\begin{equation}\label{limsup=v}
\lim_{n\to\i}\sup_{f\in \cC_n}\massE[U(f)-fyZ_T^0]=\sup_{x>0}\{u(x,U)-xy\},
\end{equation}
where $\cC_n:=\{f\in  \LL_+^0|0 \leq f \leq n\}$  is the ball of radius $n$ in the positive orthant of $\LL^0$.
Notice  the left-hand side of \eqref{limsup=v} is an increasing limit in $n$,
 so we only need to show   for each $n$ and each $f\in \cC_n,$
$$\massE[U(f)-fyZ_T^0] \leq \sup_{x>0}\{u(x,U)-xy\}.$$
To do so,
fix $f$ and define $x^*:=\massE[fZ_T^0]$.
For $x^*=0$, we have $0=\massE[fZ_T^0]$, hence $f \equiv 0$ and so $\massE[U(f)-fyZ_T^0] =U(0) \leq u(x,U)-xy+ xy$ for any $x>0$. This gives the above inequality,
so we need to consider the case of $x^*>0$.
By definition of $x^*$, we have $f \in \cC(x^*)$ and it follows that
$$\massE[U(f)-fyZ_T^0] \leq  u(x^*,U)-x^*y  \leq \sup_{x>0}\{u(x,U)-xy\},$$
which proves \eqref{limsup=v}.
  
To get \eqref{defv} from \eqref{limsup=v}, we now want to interchange supremum and expectation on the left-hand side of \eqref{limsup=v} and then let $n \to \i$. 
For each $n$, a measurable selection argument (see \cite[Theorem 18.19]{AB06}) shows that we can choose a measurable
selector $ x^*(\omega)\in \lzp$  such that 
$$\sup_{0 \leq x \leq n}\{U(x)-xyZ_T^0(\omega)\}=U(x^*(\omega))-x^*yZ_T^0(\omega).$$
With $V_n(y):=\sup_{0 \leq x \leq n}\{U(x)-xy \} \geq U(n)-ny$,
 it thus follows that we have indeed 
\begin{equation}\label{supn=Evn}
\sup_{f\in \cC_n}\massE[U(f)-fyZ_T^0]=\massE[V_n(yZ_T^0)].
\end{equation}
But $V_n$ is increasing in $n$ and dominated by $V$, so in view of \eqref{limsup=v} and \eqref{supn=Evn}, we have to show  for \eqref{defv} that $\lim_{n\to\i}\massE[V_n(yZ_T^0)]\geq \massE[V(yZ_T^0)]$.
Because of $V_n^{\pm} \to V^{\pm}$ and Fatou's lemma,
it is sufficient to show that $(V_n^{-}(xy))_{n\in \NN}$ is u.i.
Since $V_n$ is increasing in $n$, $V_n^-$ is decreasing in $n$ and  $V_1^-(yZ_T^0) \leq |U(1)|+yZ_T^0$is an integrable upper bound for $V_n^{-}(xy), n\in \NN.$

      It remains to prove that $u(x,U_c )$ and the concave envelope of $u(x,U)$ coincide on $(0,\i)$. 
  We also refer to   Part (3) from the proof of Theorem 4.1 in \cite{Rei13}.
Note first from Lemma \ref{chaV} that $U$ and $U_c$ have the same convex conjugate $V$ so that applying \eqref{sameconvexconj} for $U$ and  $U_c$ implies that  $v$ is the conjugate of both $u(x,U)$  and $u(x,U_c)$.
Thus also their biconjugates coincide. 
But applying Lemma \ref{helplemma} (iii) to $u(x,U_c )$ gives that  
$u(x,U_c )$ is equal to its biconjugate. This shows that $u(x,U_c )$ is the biconjugate of $u(x,U)$.
So applying  Lemma \ref{helplemma} (iii) to $u(x,U)$, which is upper-semicontinuous due
to Proposition \ref{upsemicont},  gives that $u(x,U_c )$ and the concave envelope of $u(x,U)$
coincide on $(0,\i).$ 
\end{proof}

\section{Appendix}
Let $f: \mathbb R \to  \bar{ \mathbb R}=\mathbb R \cup \{\pm \infty\} $ be an extended real-valued function.   
The {\it effective domain} of  $f$  is the set $\{x \in \mathbb R | f(x)<\infty       \}$ denoted by $\dom(f)$,
and its interior is denoted by  $int(\dom(f))$. The function $f$  is called {\it proper} if both $\dom(f) \neq \emptyset$ and $f(x)>-\infty$, for all $x$. 
The {\it conjugate of $f$} is the  extended real-valued function $f^*$ on $\mathbb R$ defined by
\begin{equation}
f^*(y):=\sup_{x \in \mathbb R}\{xy-f(x)\}, \quad \forall y\in \mathbb R.
\end{equation}
The {\it biconjugate } $f^{**}$  is defined by $f^{**}:=(f^*)^*$.
If $f$ is proper, lower-semicontinuous and convex, then its subdifferential $\partial f$ is the multivalued mapping defined by $\partial f(x):=\emptyset$ if $f(x)=\infty$ and 
$$\partial f(x):=\{y\in \RR | f(a)\geq f(x)+y(a-x), \forall a\in \RR \}, \mbox{if }x \in \dom(f).$$
The {\it convex envelope } $\bar f$ of $f$ is the largest convex function $\bar f \leq f.$
 
\begin{lemma}\label{helplemma}
   Suppose that $f$ is proper and lower-semicontinuous and its convex envelope $\bar f$  is proper as well. Then:
\begin{enumerate}[(i)]
\item $f$ is convex, proper and lower-semicontinuous.
\item $f$ and its convex envelope have the same conjugate.
\item  $f^{**}$is the lower-semicontinuous envelope of the convex envelope $\bar f$.
\item   Fix $x_0 \in \dom(f) $ and $y_0 \in \partial f^{**}(x_0)$. If the conjugate $f^*$ is differentiable in $y_0$ , then $f(x_0) = f^{**}(x_0)$.
\end{enumerate}

 Let $f$   in addition  be convex and extend the right and left derivative functions $f'_+$ and $f'_-$ beyond the interval $\dom(f)$ by setting both $= \infty$ for points lying to the right of $\dom(f)$  and both $= -\infty$ for points lying to the left. Then:
\begin{enumerate}
\item[(v)]  $y \in \partial f(x)$ if and only if $x \in \partial f^*(y)$.
\item[(vi)]   $y \in \partial f(x)$ if and only if $f(x)+f^*(y)=xy$.
\item[(vii)]   Let $x \in \dom( f )$. $f$ has a unique subgradient at $x$ if and only $f$ is differentiable at $x$.
\item[(viii)]   $\partial f(x)=\{y\in \RR|f'_-(x) \leq y \leq f'_+(x) \}.$
\item[(ix)]  $f'_+(z_1) \leq f'_-(x) \leq  f'_+(x)\leq  f'_-(z_2)$ when $z_1 < x <z_2.$
\end{enumerate}

 \begin{proof}
  Statements (i),
(ii) and (iii) are part of \cite[Theorem 11.1]{RW98}. Part (iv) can be found in
a similar form in Theorem 2 of   \cite{Str09}. The other properties can be found in \cite{Roc70}, see Theorem 23.5  for (v) and (vi), Theorem 25.1 for (vii), Theorem 24.1 and its discussion for (viii) and (ix).
\end{proof}

\end{lemma}

\bibliography{nonconcave}
\bibliographystyle{alpha}

\end{document}